\documentclass[12pt]{article}
\usepackage{amssymb}
\usepackage{amsmath}
\newtheorem{theorem}{Theorem}[section]
\newtheorem{lemma}{Lemma}[section]
\newtheorem{corollary}{Corollary}[section]
\newtheorem{proposition}{Proposition}[section]

\newcommand{\qed}{\hfill $\square$\vskip .2cm}

\newcommand{\sect}[1]{\section{#1}\setcounter{equation}{0}}
\newcommand{\be}{\begin{equation}}
\newcommand{\ee}{\end{equation}}

\def\v{\vskip .2cm}

\def\F{{\cal{F}}}

\def\B{{\cal{B}}}
\def\X{{\bf X}}

\begin{document}
 \title{A Martingale Proof of Dobrushin's Theorem 
 for Non-Homogeneous Markov Chains}

 \author{S. Sethuraman \ \ and \ \ S.R.S. Varadhan\\
Iowa State University\ \ and \ \ New York University}

 \thispagestyle{empty}
 \maketitle
 \abstract{In 1956, Dobrushin proved a definitive
 central limit theorem for
   non-homogeneous Markov chains.  In this note, a shorter and
   different proof elucidating
   more the assumptions is given through martingale approximation.

 }  

 \v
  \thanks{

\noindent Partially supported by NSF/DMS-0071504 and NSF/DMS-0104343. \\
E-mail: sethuram@iastate.edu and varadhan@cims.nyu.edu\\[.15cm]
 \\[.15cm]
 {\sl Abbreviated title}:  Martingale Proof of Dobrushin's Theorem
 \\[.15cm]
 {\sl AMS (2000) subject classifications}: 
 Primary 60J10; secondary 60F05.}

\section{Introduction and Results}
Nearly fifty years ago,
R. Dobrushin proved in his thesis \cite{Dobrushin} 
a definitive central limit
theorem (CLT) for  Markov chains in discrete time that are not necessarily
homogeneous in time.  Previously, Markov, 
Bernstein, Sapagov, and Linnik, among others, had 
considered the central limit question 
under various sufficient conditions.  Roughly, the progression of
results relaxed the state space structure from $2$ states to an
arbitrary set of states, and also the level of asymptotic degeneracy
allowed  for the transition probabilities of the chain.

After Dobrushin's
work, some refinements and extensions of his CLT, 
some of which under more stringent
assumptions, were proved by Statulevicius \cite{Statulevicius} and Sarymsakov 
\cite{Sarmysakov}. See also Hanen \cite{H} in this regard.  
A corresponding
invariance principle was also
proved by Gudinas \cite{Gudynas}. More general references on non-homogeneous
Markov processes can be found in 
Isaacson and Madsen \cite{Isaacson}, Iosifescu \cite{I}, Iosifescu and Theodorescu 
\cite{Iosifescu}, and Winkler \cite{Winkler}.


We now define what is meant by ``degeneracy.''  
Although there are many 
measures of ``degeneracy,'' the measure which turns out to be most
useful to work with is that in terms of the contraction coefficient.
This coefficient has appeared in early results concerning
Markov chains, however, in his thesis, Dobrushin popularized its use,
and developed many of its important properties. [See Seneta
\cite{Seneta} for some history.]  

Let $\pi = \pi(x, dy)$ be  a Markov transition probability on $(\X, {\mathcal B}(\X))$. Define the contraction coefficient $\delta(\pi)$ of $\pi$ as
\begin{align*}
\delta(\pi)&= \sup_{x_1,x_2\in \X, A\in {\mathcal B}}|\pi (x_1, A)-\pi(x_2,A)|\\
            &={1\over 2} \sup_{f:|f|\le 1}|\int_\X f(y)[\pi(x_1, dy)-\pi(x_2, dy)]|
\end{align*}
Also, define the related coefficient
$\alpha(\pi)=1-\delta(\pi)$.

Clearly, $0\leq \delta(\pi)\leq 1$, and $\delta(\pi)=0$ if and only if  $\pi(x, dy)$ is independent of $x$. It makes sense to call $\pi$ ``non-degenerate'' if $0\le \delta(\pi)<1$. We use the standard convention
and denote by $\mu\pi$ and $\pi u$ the transformations induced by $\pi$ on  countably additive signed  measures and bounded measurable functions respectively,
$$
(\mu\pi)(A)=\int \pi(x,A)\,\mu(dx)\qquad{\rm and }\qquad (\pi u)(x)=\int u(y)\,\pi(x,dy)
$$

It is easy to see that $\delta (\pi)$ has the following properties.
$$
\delta(\pi)=\sup_{x_1,x_2\atop  u\in {\mathcal U}}|(\pi u)(x_1)-(\pi u)(x_2)|
$$
where ${\mathcal U}=\{u:\sup_{y_1,y_2}|u(y_1)-u(y_2)|\le 1\}$.
It is the operator norm of $\pi$  with respect to the  Banach (semi-)
norm ${\rm Osc}(u)=\sup_{x_1, x_2}
|u(x_1)-u(x_2)|$, namely the oscillation of $u$. In particular, for any transition probabilities
$\pi_1, \pi_2$ we have
\begin{equation}\label{product estimate}
\delta(\pi_1\,\pi_2)\le \delta(\pi_1)\,\delta(\pi_2).
\end{equation}
 If $\mu$ is a signed measure with $\mu(\X)=0$,
$$
\|\mu\|_{\rm Var}=\sup_A|\mu(A)|=\frac{1}{2}\sup_{\|u\|_{L^\infty}\le 1}|\int f(x) \mu(dx)|=\sup_{u\in{\mathcal U}}|\int u(x)\,\mu(dx)|.
$$
Therefore, by duality, for any two probability measures $\lambda$ and $\mu$ on $\X$,
\begin{equation}\label{dual estimate}
\|(\lambda-\mu)\pi\|_{\rm Var}\le \delta(\pi)\|\lambda-\mu\|_{\rm Var}.
\end{equation}

By a non-homogeneous Markov chain of length $n$  on state space $(\X, {\mathcal B}(\X))$ 
corresponding to transition operators
$\{\pi_{i,i+1}=\pi_{i,i+1} (x, dy):  1\le i\le n-1\}$ we mean the Markov process $P$
on the product space 
$(\X^n, \ \B(\X^n))$, 
$$P[X_{i+1}\in A|X_{i} = x] = \pi_{i,i+1}(x, A),$$
where $\{X_i: 1\le i\le n\}$ are the canonical projections. In particular, under the initial distribution $X_1\sim\mu$, the distribution at
time $k\geq 1$ is $\mu \pi_{1,2}\pi_{2,3}\cdots \pi_{k-1,k}$. 
For $i<j$ we will define
$$
\pi_{i,j} = \pi_{i,i+1}\pi _{i+1,i+2}\cdots \pi_{j-1,j}
$$ 
We denote  by $E[Z]$ and $V(Z)$ the
expectation and variance of the random variable $Z$ with respect to $P$.

Consider now a non-homogeneous Markov chain on $\X$ with respect
to transition operators $\{\pi_{i,i+1}: 1\le i\le n-1\}$.  
The following comparison of
marginal distributions at time $n$ starting from different initial
conditions is an easy consequence of (\ref{product estimate}) and (\ref{dual estimate}).

\begin{eqnarray}
\|\lambda \pi_{1,n} - \mu \pi_{1,n}\|_{\rm Var} & \leq &
\|\mu -
\nu\|_{\rm Var} \ \delta(\pi_{1,n})\nonumber\\
&\leq& \|\mu -\nu\|_{\rm Var}\ \prod_{i=1}^{n-1} \delta(\pi_{i,i+1}).\label{nh_var_ineq}
\end{eqnarray}

Dobrushin's 
theorem concerns the fluctuations of an array of
non-homogeneous Markov chains.
For each $n\geq 1$, let $\{X_i^{(n)}: 1\leq i\leq n\}$ 
be  $n$ observations
of a non-homogeneous Markov chain on 
$\X$ with transition matrices 
$\{\pi^{(n)}_{i,i+1} = \pi^{(n)}_{i,i+1}(x, dy): 1\leq i\leq n-1\}$.  Let also 
$$\alpha_n = \min_{1\leq i\leq n-1} \alpha\big(\pi^{(n)}_{i,i+1}\big).$$
In addition, let $\{f^{(n)}_i: 1\leq i\leq n\}$ be real valued functions
on $\X$.  Define, for $n\geq 1$, the sum
$$S_{n} = \sum_{i=1}^n f^{(n)}_l(X_i^{(n)}).$$

\begin{theorem}
\label{mt_dob}
Suppose that for some finite constants $C_n$, 
$$\sup_x\sup_{1\le i\le n} |f^{(n)}_i(x)| \leq C_n.$$   Then, if
\be
\label{alpha_cond}
\lim_{n\rightarrow \infty}
C^2_n\alpha^{-3}_n \bigg[\sum_{i=1}^n
V\big(f^{(n)}_i(X^{(n)}_i)\big)\bigg]^{-1} = 0,
\ee
we have, regardless of the initial distribution, that
\be
\label{convergence}
\frac{S_n - E[S_n]}{\sqrt{V(S_n)}} \ \Rightarrow \ {\rm N}(0,1).
\ee
In general, the result is not true if condition
(\ref{alpha_cond}) is not met.
\end{theorem}

In \cite{Dobrushin}, Dobrushin also states the direct corollary which
simplifies some of the assumptions.
\begin{corollary}
\label{maincor}
When the functions are uniformly bounded, i.e.  $\sup_n C_n=C<\infty$
and the variances
are bounded below, i.e.
$V(f^{(n)}_i(X^{(n)}_i))\geq c>0$, for all $1\leq i\leq n$ and
$n\geq 1$, then we have the
convergence (\ref{convergence}) provided
\be
\label{corcond}\lim_{n\rightarrow \infty}
n^{1/3}\alpha_n =\infty.
\ee
\end{corollary}

We remark that 
in \cite{Dobrushin} (e.g. Theorems 3, 8) there are also results 
where the boundedness
condition on $f^{(n)}_i$ is replaced by integrability conditions.  As
these results follow from truncation methods and Theorem \ref{mt_dob}
for bounded variables, we only consider
Dobrushin's theorem in the bounded case.  

Also, for the ease of the reader, and to be complete, we will
discuss in the next section an example, given in \cite{Dobrushin} and
due to Dobrushin and Bernstein,
of how
the weak convergence (\ref{convergence}) may fail when the condition
(\ref{alpha_cond}) is not satisfied.

We now consider Dobrushin's methods.  
The techniques used in \cite{Dobrushin} to prove the above
results fall under the general heading of the
``blocking method.''  The condition (\ref{alpha_cond}) ensures that
well-separated blocks of observations may 
be approximated by independent versions with small error.  
Indeed, in many remarkable steps, Dobrushin
exploits the Markov
property and several contraction coefficient properties, which he
himself derives, to deduce error bounds sufficient to apply CLT's for
independent variables.
However, in \cite{Dobrushin},
it is difficult to see,
even at the technical level, 
why condition (\ref{alpha_cond}) is natural.

The aim of this note is to provide a different, shorter proof of
Theorem \ref{mt_dob} which explains more why
condition (\ref{alpha_cond}) appears in the result.  The methods are
through martingale approximations and martingale CLT's which perhaps were not as codified
in the early 1950's
as they are today.  These methods go back at least to Gordin \cite{Gordin}
in the context of homogeneous processes,
and has been used by others in other ``related'' situations
(e.g. Kifer \cite{Kifer}; see also Pinsky \cite{Pinsky}).  There are three main ingredients 
in this approximation with respect to 
the non-homogeneous setting of Theorem \ref{mt_dob}, (1) negligibility estimates for individual
components, (2) a law of large numbers for conditional variances, 
and (3) lower
bounds for the variance $V(S_n)$.  Negligibility bounds and a
LLN are
well known requirements for martingale CLT's 
(cf. Hall-Heyde \cite[ch. 3]{Hall-Heyde}), and in fact, as will be
seen, the sufficiency of condition (\ref{alpha_cond}) is transparent
in the proofs of
these two components (Lemma \ref{negligibility}, and Lemmas \ref{osc}
and \ref{nh_osc}).
The variance lower bounds which we will use 
were as well
derived
by Dobrushin in his proof.  However, using some martingale properties, we give
a more direct argument for a better estimate.

We note also, with this martingale approximation, that
an invariance principle for the partial sums holds through
standard martingale propositions, Hall-Heyde \cite{Hall-Heyde}, among other results.  In fact, from 
the martingale invariance principle,
it should be 
possible to derive Gudynas's
theorems
\cite{Gudynas} although this is not done here.

We now explain the structure of the article.  In section 2, we give the
Bernstein-Dobrushin example of a Markov chain with anomalous
behavior.  In section 3, we discuss needed properties of the
contraction coefficient.  In section 4, we 
state the martingale CLT that will be utilized, and, as a preview of
the
non-homogeneous chain proof, we quickly reprise the argument with
respect to homogeneous
chains.  In section 5, we prove Theorem \ref{mt_dob} with martingale
approximation assuming a lower bound on the variance
$V(S_n)$.  And last, in section 6, we prove this variance
estimate.  

\sect{Anomalous Example}

Here, we summarize the example in Dobrushin's thesis,
attributed to Bernstein, which shows that condition
(\ref{alpha_cond}) is sharp.
\vskip .2cm

\noindent {\bf Example 2.1}\ 
Let $\X = \{1,2\}$, and consider 
the $2\times 2$ transition matrices on $\X$,
$$Q(p)=\left(\begin{array}{cc}1-p&p\\p&1-p\end{array}\right)$$
The contraction coefficient $\delta(Q(p))$ of $Q(p)$ is $|1-2p|$. Note that
$\delta(Q(p))=\delta(Q(1-p))$.
The invariant  measures for all the $Q(p) $
are the same $p(1)=p(2)=\frac{1}{2}$.  
We will be looking at $Q(p)$ for $p$ close to $0$ or $1$ and the special case of $p=\frac{1}{2}$.
However,  when $p$ is small, the
homogeneous chains behave  very differently under $Q(p)$  and $Q(1-p)$.  More
specifically, when $p$ is small
there are very few  switches between the two states whereas when $1-p$ is small it switches most of the time. In fact, this behavior can be made more precise (see Dobrushin \cite{D}, 
Hanen \cite{H}, or from direct computation).  
Let $T_n= \sum_{i=1}^n {\bf 1}_{\{1\}}(X_i)$
count the number of visits to state $1$, say, in 
$n$ steps.  
\vskip .2cm
{\it Case A.}  Consider the homogeneous chain under $Q(p)$ with $p=\frac{1}{n}$ and  initial distribution $p(1)=p(2)=\frac{1}{2}$.  Then,
\be
\label{A}
\frac{T_n}{n}\ \Rightarrow \ G \ \ \ {\rm and \ \ \ }
\lim_{n\to\infty }n^{-2} \ V(T_n) \ = \ V_A<\infty.
\ee
where $G$ is a proper distribution supported on $[0,1]$.

\vskip .2cm
{\it Case B.}  Consider the homogeneous chain run
under  $Q(p)$ with $p=1-\frac{1}{n}$ and initial distribution  $p(1)=p(2)=\frac{1}{2}$. Then,
\be
\label{B}
T_n -\frac{n+1}{2} \ \Rightarrow \ F
\ \ \ {\rm and \ \ \ }
\lim_{n\rightarrow \infty} V(T_n) \ = \ V_B<\infty.
\ee
where $F$ is a proper distribution function.

\vskip .1cm

Let a sequence $\alpha_n\to 0$ with $\alpha_n\ge n^{-{1\over 3}}$ be given . To construct the anomalous Markov chain, it will be helpful to split
the time horizon $[1,2,\ldots,n ]$ into roughly $n\alpha_n$ blocks of size $\alpha_n^{-1}$. 
We interpose a $Q(\frac{1}{2})$ between any two blocks that has the effect of making the blocks independent of each other.  More precisely
let
$k^{(n)}_i = i[\alpha_n^{-1}]$ for $1\leq i \leq m_n$
where $m_n=
[n/[\alpha_n^{-1}]]$.  Also, define $k^{(n)}_0=0$, and
$k^{(n)}_{m_n+1}=n$.

Define now, for $1\leq i\leq n$,
$$\pi_{i,i+1}^{(n)} = \left\{ \begin{array}{rl} Q(\alpha_n)&{\rm \ for \
      }i=1,2,\ldots,k^{(n)}_1 - 1\\
Q(\frac{1}{2})&{\rm \ for \ }i=k^{(n)}_1,k^{(n)}_2,\ldots,k^{(n)}_{m_n}\\
Q(1-\alpha_n)&{\rm \ for \ all \ other \ }i.\end{array}\right.
$$

Consider the non-homogeneous chain with respect to $\{\pi^{(n)}_{i,i+1}: 1\leq i\leq
n-1\}$ starting from equilibrium $p(0)=p(1)=\frac{1}{2}$.  
From the definition of the chain, one observes, as $Q(\frac{1}{2})$ does not
distinguish between states, that the process in time horizons $\{(k^{(n)}_i +1,
k^{(n)}_{i+1} ): 0\leq i\leq m_n\}$ are mutually independent.  For the first
time segment $1$ to $k^{(n)}_1$, the chain is in regime $A$, while
for the other segments, the chain is in case $B$.

Once again, let us concentrate on 
the number of visits to 
state $1$.  Denote by $T^{(n)} = \sum_{i=1}^n {\bf 1}_{\{1\}}(X^{(n)}_i)$ 
and $T^{(n)}(k,l)= \sum_{i=k}^{l} {\bf 1}_{\{1\}}(X^{(n)}_i)$ 
the counts in the
first $n$ steps and in steps $k$ to $l$ respectively.
It follows from the discussion of independence above that
$$T^{(n)} \ = \ \sum_{i=0}^{m_n}  T^{(n)}(k^{(n)}_i+1, k^{(n)}_{i+1})$$
is the sum of independent sub-counts where, additionally, the
sub-counts for $1\leq i \leq m_n-1$ are identically distributed, the
last sub-count perhaps being shorter.
Also, as the initial distribution  is invariant, we have $V({\bf 1}_{\{1\}}(X^{(n)}_i))=1/4$
for all $i$ and $n$.  Then, in the notation of Corollary \ref{maincor},
$C=1$ and $c=1/4$.

From (\ref{A}), we have that
$$V(T^{(n)}(1,k^{(n)}_1)) \ \sim \ \alpha^{-2}_n V_A \ \ {\rm
  as \ }n\uparrow \infty.$$
Also, from (\ref{B}) and independence of $m_n$ sub-counts, we have that
$$V(T^{(n)}(k^{(n)}_1+1,n)) \ \sim\ n\alpha_n V_B \ \
{\rm as \ }n\uparrow \infty.$$
From these calculations, we see if $n^{1/3}\alpha_n\rightarrow
\infty$, then $\alpha^{-2}_n << n \alpha_n$, and so the major
contribution to $T^{(n)}$ is from $T^{(n)}(k^{(n)}_1+1,n)$.
However, since this last count is (virtually) 
the sum of $m_n$ i.i.d. sub-counts,
we have that $T^{(n)}$, properly normalized, converges to
${\rm N}(0,1)$, as predicted by Dobrushin's theorem
\ref{mt_dob}.

On the other hand, if $\alpha_n = n^{-1/3}$, we have $\alpha^{-2}_n =
n\alpha_n$, and count $T^{(n)}(1,k^{(n)}_1)$, 
independent of $T^{(n)}(k^{(n)}_1,n)$, also contributes to the
sum $T^{(n)}$.  After appropriate scaling, then, $T^{(n)}$
approaches the convolution of a non-trivial 
non-normal  distribution and a normal distribution,  and therefore  is certainly not
Gaussian.  


\sect{Martingale CLT}

The central limit theorem for martingale differences is by now a standard
tool.  We quote the following (strong) 
form of the result implied
by Corollary 3.1 in 
Hall and Heyde \cite{Hall-Heyde}.
\begin{proposition}
\label{martprop}
For each $n\ge 1$, let $\{(W^{(n)}_i,\F^{(n)}_i): 0\leq i\leq n\}$
be a martingale  relative to the nested family 
$\F^{(n)}_i\subset \F^{(n)}_{i+1}$ with  $W^{(n)}_0=0$. Let 
$\xi^{(n)}_i=W^{(n)}_i-W^{(n)}_{i-1}$ be their differences.
Suppose that
$$\begin{array}{rll}
\max_{1\leq i\leq n} \|\xi^{(n)}_i\|_{L^\infty} & \rightarrow\  0 \ &\
{\rm and }\\
\sum_{i=1}^n E[(\xi^{(n)}_i)^2|\F^{(n)}_{i-1}]&\rightarrow \ 1 \ &\ {\rm in\ }L^2.
\end{array}$$
Then, 
$$W^{(n)}_n \ \Rightarrow \ {\rm N}(0,1).$$
\end{proposition}
Note that 
the first and second limit conditions are the 
negligibility assumption on the
sequence, and
law of large numbers for conditional
variances mentioned in the introduction.

We now sketch a proof of  Corollary \ref{maincor}  in the 
case  of a homogeneous Markov chain on a finite state space. Assume that we have a Markov chain with transition probability $P$ on a finite state space $\X$. If $\delta(P)<1$, and $f: { \X}\to R$ is a function with mean $0$  with respect to the invariant distribution $\pi$ on $ \X$, it is in the range of $I-P$ and the equation $(I-P)u=f$ has a solution.  The following argument is implicit  in Gordin \cite{Gordin}, and also explicitly used  in Kipnis and Varadhan \cite{KV}.

Using the relation $E[u(X_{j+1})|{\cal F}_{j}]=(Pu)(X_{j})$, it is easy to check that
\begin{equation*}
f(X_j)= u(X_j)-E[u(X_{j+1})|{\cal F}_j]=u(X_j)-u(X_{j+1})+\xi_{j+1}
\end{equation*}
where
\begin{equation*}
\xi_j=u(X_{j})-E[u(X_{j})|{\cal F}_{j-1}]
\end{equation*}
is a martingale difference. Then,
\begin{equation*}
\sum_{j=0}^{n-1} f(X_j)=u(X_0)-u(X_{n})+\sum_{j=1}^n \xi_j .
\end{equation*}
If we define 
$$
E[\xi_{j+1}^2|{\cal F}_{j}]=q(X_{j})
$$
We will apply the martingale CLT (Proposition \ref{martprop}) to the array
formed from $W^{(n)}_i = M_i/\sqrt{n}$ with differences $\xi^{(n)}_i =
\xi_i/\sqrt{n}$.  As the differences are uniformly bounded,
$\|\xi^{(n)}_i\|_{L^\infty} \leq 2\|u\|_{L^\infty}/\sqrt{n}$, the
first condition of Proposition \ref{martprop} is satisfied.  The second
follows from the following computation.
From the Markov property, 
$$
\sum_{i=1}^n E[(\xi^{(n)}_i)^2 |\F^{(n)}_{i-1}] \ =\ 
\frac{1}{n}\sum_{i=1}^n E[\xi^2_i|X_{i-1}]
\ =\ \frac{1}{n}\sum_{i=0}^{n-1} q(X_i).$$
So, by
the ergodic theorem, the last expression converges almost surely to $V_0
= E_\pi[q(X_0)]<\infty$. It is not difficult to see that $V_0>0$.

Therefore, $V(M_n) \sim nV_0$ and $(nV_0)^{-1/2}M_n \Rightarrow {\rm
  N}(0,1)$ by Proposition \ref{martprop}. Since 
the difference 
$$n^{-{1\over 2}}|M_n-\sum_{j=0}^{n-1} f(X_j)|\le 
n^{-1/2}\|u(X_n) -u(X_0)\|_{L^\infty} \to 0
$$  we have
$V(S_n)\sim nV_0$ and
$S_n/\sqrt{V(S_n)}\Rightarrow {\rm N}(0,1)$ also.
\qed

\sect{Proof of Theorem \ref{mt_dob}}
We give here a short proof for Theorem \ref{mt_dob}
through 
martingale approximation, 
illustrated for homogeneous chains in the previous section.
Consider the 
non-homogeneous setting of Theorem \ref{mt_dob}.  To follow the
homogeneous argument, we 
will need to find the non-homogeneous analogue of
the resolvent function ``$u = (I-P)^{-1}f$.'' 
To simplify notation, we will assume throughout
that the functions $\{f^{(n)}_i\}$ are mean-zero, 
$E[f^{(n)}_i(X^{(n)}_i)]=0$ for $1\leq i\leq n$ and $n\geq 1$.  
Define
\begin{equation*}
Z^{(n)}_k = \sum_{i=k}^n E[f^{(n)}_i(X^{(n)}_i)|X^{(n)}_k] 
\end{equation*}
where
\begin{equation}
\label{gen_resolvent}
Z^{(n)}_k =
\begin{cases}
f^{(n)}_k(X^{(n)}_k) + \sum_{i=k+1}^n
E[f^{(n)}_i(X^{(n)}_i)|X^{(n)}_k] & {\rm for}\quad  1\le k\le n-1\\
f^{(n)}_k(X^{(n)}_n) & {\rm for}\quad    k=n.
\end{cases}
\end{equation}
\vskip .1cm

\noindent {\bf Remark 4.1} \ Before going further, we 
remark that indeed sequence $\{Z^{(n)}_k\}$  can be thought of as
a generalization of
the resolvent sequence $\{u(X_k)\}$ used in the case of a
homogeneous chain.
When the array $\{X^{(n)}_i\}$ is formed from the sequence $\{X_i\}$, 
$f^{(n)}_i = f$ for all $i$ and $n$, and the chain
is homogeneous, $P_n = P$ for all $n$, then indeed $Z^{(n)}_k$ reduces to
$Z^{(n)}_k = 
f(X_k) + \sum_{i= 1}^{n-k} (P^i f)(X_k)$
which approximates
$\sum_{i=0}^\infty (P^if)(X_k) 
 = [(I-P)^{-1}f](X_k)  =   u(X_k)$.
See also p. 145-6 Varadhan \cite{Varadhan} for other uses of $\{Z^{(n)}_k\}$.

Now, let us return to the full non-homogeneous setting of Theorem
\ref{mt_dob}.  By 
rearranging terms in (\ref{gen_resolvent}), we obtain for $1\leq k\leq n-1$
\begin{eqnarray}
\label{rearrange}
f^{(n)}_k(X^{(n)}_k) &=& Z^{(n)}_k -E[Z^{(n)}_{k+1}|X^{(n)}_k]\\
&=&\big[Z^{(n)}_k -E[Z^{(n)}_k|X^{(n)}_{k-1}]\big] + \big[E[Z^{(n)}_k|X^{(n)}_{k-1}] -
E[Z^{(n)}_{k+1}|X^{(n)}_k]\big].\nonumber
\end{eqnarray}
Then, we have the decomposition,
\begin{eqnarray}
\label{Sn_decomp}
S_n &=& \sum_{k=1}^n f^{(n)}_k(X^{(n)}_k) \nonumber\\
&=& \sum_{k=2}^n [Z^{(n)}_k - E[Z^{(n)}_k|X^{(n)}_{k-1}]]
+Z^{(n)}_1
\end{eqnarray}
and so in particular $V(S_n) = \sum_{k=2}^n V(Z^{(n)}_k -
E[Z^{(n)}_k|X^{(n)}_{k-1}]) + V(Z_1^{(n)})$.
Let us now define the differences
\be
\label{xi_decomp}
\xi^{(n)}_k = \frac{1}{\sqrt{V(S_n)}} 
\big[Z^{(n)}_k - E[Z^{(n)}_k|X^{(n)}_{k-1}]\big]\ee
and the martingale $M^{(n)}_k = \sum_{l=2}^k \xi^{(n)}_l$ with respect
to $\F^{(n)}_k = \sigma\{X^{(n)}_l: 1\leq l\leq k\}$ for $n\geq 1$.
The plan to obtain Theorem \ref{mt_dob} will now be to approximate
$S_n$ by $M^{(n)}_n$ and use 
Proposition \ref{martprop}.
Condition (\ref{alpha_cond}) will be a natural sufficent condition for ``negligibility'' (Lemma
\ref{negligibility}) and ``LLN'' (Lemmas \ref{osc} and \ref{nh_osc})
with regard to
Proposition \ref{martprop}.


\begin{lemma}
\label{Z_bound}
For $1\leq i\leq j\leq n$, we have the bound
$$\|E[f^{(n)}_j(X^{(n)}_j)|X^{(n)}_i]\|_{L^\infty}\leq 2C_n(1-\alpha_n)^{j-i}.$$
Hence, for $1\leq k\leq n$, 
$$\|Z^{(n)}_k\|_{L^{\infty}} \ \leq \
2C_n\alpha_n^{-1}.$$
\end{lemma}

{\it Proof.}  Since $\|f_j^{(n)}\|_{L_\infty}\le C_n$ its oscillation
${\rm Osc}(f_j^{(n)})\le 2C_n$.  From (\ref{nh_var_ineq}),
$$
{\rm Osc}(\pi_{i,j}f_j^{(n)})\ \le \ 2C_n \delta (\pi_{i,j})\ \le \ 2C_n (1-\alpha_n)^{j-i}
$$
Because $E[(\pi_{i,j}f_j^{(n)})(X_i)]=E[f_j^{(n)}(X_j)]=0$, 
$$
\|\pi_{i,j}f_j^{(n)}\|_{L_\infty}\ \le \ {\rm
  Osc}(\pi_{i,j}f_j^{(n)})\ \le \ 2C_n (1-\alpha_n)^{j-i}
$$
The second estimate now follows from this estimate.
Indeed,
$$
|Z^{(n)}_k| \ =\  |\sum_{i=k}^n [E[f^{(n)}_i(X^{(n)}_i)|X^{(n)}_k]|
\ \leq\  2C_n\sum_{i=k}^n (1-\alpha_n)^{k-i}
\ \leq\ 
2C_n\alpha_n^{-1}. \ \ \ \ \ \ \ \ \ \ \ \ \ \ \ \square
$$

We now state a lower bound for the variance which will be proved in the next
section using martingale ideas.
We remark in
\cite{Dobrushin} that actually the bound, $V(S_n) \geq
(\alpha_n/8)\sum_{i=1}^n V(f^{(n)}_i(X^{(n)}_i))$, is found by
different methods
(see also section 1.2.2 \cite{Iosifescu}). 
\begin{proposition}
\label{lowerbound}
For $n\geq 1$,
\be
\label{lb}
V(S_n) \ \geq\
\frac{\alpha_n}{4}\sum_{i=1}^n V\big(f^{(n)}_i(X^{(n)}_i)\big).
\ee
\end{proposition}

The next estimate shows that
the asymptotics of $S_n/\sqrt{V(S_n)}$ 
depend only on
the martingale approximant $M^{(n)}_n$,
and
that the differences $\xi^{(n)}_k$ are negligible.  


\begin{lemma}
\label{negligibility}
Under condition (\ref{alpha_cond}),
we have that
$$\lim_{n\rightarrow \infty}
\sup_{1\leq k\leq n}\frac{\|Z^{(n)}_k\|_{L^\infty}}{\sqrt{V(S_n)}} = 0.$$
\end{lemma}

{\it Proof.}
By Lemma \ref{Z_bound} and Proposition \ref{lowerbound}, 
$$\frac{\|Z^{(n)}_k\|_{L^\infty}}{\sqrt{V(S_n)}} 
\leq \frac{4C_n}{\big(\alpha^3_n\sum_{i=1}^n
  V\big(f^{(n)}_i(X^{(n)}_i)\big)\big)^{1/2}}.$$ 
The lemma follows now from (\ref{alpha_cond}).  
\qed

We now turn to showing the LLN part of Proposition 
\ref{martprop} for the array$\{M^{(n)}_k\}$.  
\begin{lemma}
\label{osc}
Let $\{Y^{(n)}_l: 1\leq l\leq n\}$ and
$\{\F^{(n)}_l:1\leq l\leq n\}$, for $n\geq 1$,
be respectively
an array of non-negative variables and $\sigma$-fields such
that $\sigma\{Y^{(n)}_1,\ldots,Y^{(n)}_l\}\subset
\F^{(n)}_l$.
Suppose 
that
$$\lim_{n\rightarrow \infty}E\big[\sum_{l=1}^nY^{(n)}_l\big]=1 {\rm \ \ \ and \ \ \ }
\sup_{1\leq i\leq n} \|Y^{(n)}_i\|_{L^\infty} \leq \epsilon_n $$
where
$\lim_{n\rightarrow\infty}\epsilon_n =0$.
In addition, assume
$$\lim_{n\rightarrow \infty}\sup_{1\leq l\leq n-1} {\rm Osc} \
E\big[\sum_{j=l+1}^n Y^{(n)}_j|\F^{(n)}_{l}\big] = 0.$$
Then, 
$$\lim_{n\rightarrow \infty}\sum_{l=1}^nY^{(n)}_l = 1 \ \ \ {\rm in \ }L^2.$$
\end{lemma}

{\it Proof.}
Write
$$
E\big[(\sum_{l=1}^n Y^{(n)}_l)^2\big] = \sum_{l=1}^n E\big[(Y^{(n)}_l)^2\big] 
+ 2\sum_{l=1}^{n-1} E\big[Y^{(n)}_l(\sum_{j=l+1}^n Y^{(n)}_j)\big].$$

The first sum on the right-hand side is bounded as follows.
From non-negativity,
$$
\sum_{l=1}^n E\big[(Y^{(n)}_l)^2\big] \ \leq\  \epsilon_n\sum_{l=1}^n
E\big[Y^{(n)}_l\big]
\ =\ \epsilon_n \cdot (1+o(1)) \ \rightarrow 0 {\rm \ \ as \ \ \ \ }n\uparrow \infty.
$$

Consider now the second sum.  Write
$$
\sum_{l=1}^{n-1} E\big[Y^{(n)}_l\big(\sum_{j=l+1}^n Y^{(n)}_j\big)\big]
=
\sum_{l=1}^{n-1} E\big[Y^{(n)}_l E\big[\sum_{j=l+1}^n Y^{(n)}_j|\F^{(n)}_l\big]\big].
$$
From the oscillation assumption, we have that
$$\sup_{1\leq l\leq n-1}\sup_\omega
\big|E\big[\sum_{j=l+1}^n Y^{(n)}_j|\F^{(n)}_l\big](\omega) -E\big[\sum_{j=l+1}^n
Y^{(n)}_j\big]\big| = o(1).$$
Therefore,
\begin{eqnarray*}
2\sum_{l=1}^{n-1} E\big[Y^{(n)}_l\big(\sum_{j=l+1}^n Y^{(n)}_j\big)\big]
&=& 2\sum_{l=1}^{n-1} E\big[Y^{(n)}_l]E[\sum_{j=l+1}^n
Y^{(n)}_j\big] + o(1)\cdot \sum_{l=1}^{n-1} E\big[Y^{(n)}_l\big]\\
&=& \big(\sum_{l=1}^{n} E\big[Y^{(n)}_l\big]\big)^2 -\sum_{l=1}^n E\big[(Y^{(n)}_l)^2\big] +
o(1)\\
&=& 1 + o(1).
\end{eqnarray*}

Putting together these statements, we obtain the lemma. \qed

To apply this result 
to our situation, we will need the following oscillation
estimate.

\begin{lemma}
\label{nh_osc}
Let $v^{(n)}_l =
E\big[(\xi^{(n)}_{l})^2|X^{(n)}_{l-1}\big]$ and
$\F^{(n)}_j=\sigma\{
X^{(n)}_1,\ldots, X^{(n)}_j\}$ for $2\leq l\leq n$
and $1\leq j\leq n$.  
Then, under condition (\ref{alpha_cond}), we have
$$\sup_{2\leq l\leq n-1}{\rm Osc} \ E\big[\sum_{j=l+1}^n 
v^{(n)}_j|\F^{(n)}_l\big] = o(1).$$
\end{lemma}

{\it Proof.} From the martingale and Markov property,
we have$E[\xi^{(n)}_r\xi^{(n)}_s|X^{(n)}_u] = 0$ for $r>s\geq u$.
Then,
\begin{eqnarray*}
E\big[\sum_{j=l+1}^nv^{(n)}_j |\F^{(n)}_l\big]
&=& E\big[\ \sum_{j=l+1}^n \big(\xi^{(n)}_{j}\big)^2\ |X^{(n)}_l\big]\\
&=& E\big[\ \bigg(\sum_{j=l+1}^n \xi^{(n)}_{j}\bigg)^2\ |X^{(n)}_l\big]\\
&=& V(S_n)^{-1}E\big[\ \big(\sum_{j=l+1}^n
f^{(n)}_{j}(X^{(n)}_{j}) - E\big[Z^{(n)}_{l+1}|X^{(n)}_l\big]\big)^2\ 
|X^{(n)}_l\big] \\
&=& (1+o(1))V(S_n)^{-1}E\big[\ \big(\sum_{j=l+1}^{n}
f^{(n)}_j(X^{(n)}_j)\big)^2\ |X^{(n)}_l\big] + o(1)
\end{eqnarray*}
where we rewrite $\xi^{(n)}_{j}$ with
(\ref{xi_decomp}) in the third line, and use Lemma \ref{Z_bound} in the
last line.

Therefore, let us consider oscillations of 
\be
\label{osc1}
E\big[\big(\sum_{j=l+1}^{n} f^{(n)}_j(X^{(n)}_j)\big)^2|X^{(n)}_l\big]
= \sum_{j,m=l+1}^{n}
E\big[f^{(n)}_j(X^{(n)}_j)f^{(n)}_m(X^{(n)}_m)|X^{(n)}_l\big].
\ee
From Lemma \ref{Z_bound}, we have the bound, for $j\leq m$,
\begin{eqnarray*}
\|E\big[f^{(n)}_j(X^{(n)}_j)f^{(n)}_m(X^{(n)}_m)|X^{(n)}_l\big]\|_{L^\infty}
&=&
\|E\big[f^{(n)}_j(X^{(n)}_j)E\big[f^{(n)}_m(X^{(n)}_m)|X^{(n)}_j\big]|X^{(n)}_l\big]\|_{L^\infty}\\
&\leq&4C^2_n(1-\alpha_n)^{j-l}(1-\alpha_n)^{m-j}.
\end{eqnarray*}
Therefore, the oscillations of (\ref{osc1}) are bounded by
$16C_n^2\alpha_n^{-2}$
uniformly in $l$.  Hence, using Proposition \ref{lowerbound},
we obtain
$$\sup_{2\leq l\leq n-1}{\rm Osc} \ E\big[\sum_{j=l+1}^nv^{(n)}_j
|\F^{(n)}_l\big] \ \leq \ (16)(4)C_n^2\bigg[\alpha^3_n\sum_{j=1}^nV\big(f^{(n)}_j(X^{(n)}_j)\big)\bigg]^{-1}$$
which is $o(1)$ by (\ref{alpha_cond}). \qed

\vskip .2cm
{\it Proof of Theorem \ref{mt_dob}.}
From Lemma \ref{negligibility}, we need only show that 
$M^{(n)}_n/\sqrt{V(S_n)}
  \Rightarrow {\rm N}(0,1)$.  This will follow from 
martingale
  convergence (Proposition \ref{martprop})
  as soon as we show (1) $\sup_{2\leq k\leq n}\|\xi_k^{(n)}\|_{L^\infty}\rightarrow
  0$ and (2) $\sum_{k=2}^n E[(\xi^{(n)}_k)^2|\F^{(n)}_{k-1}]\rightarrow 1$.
  However, (1) follows from the negligibility estimate
Lemma \ref{negligibility}, and (2) from LLN Lemmas \ref{osc}
  and \ref{nh_osc} since ``negligibility'' (1) holds and $\sum_{k=2}^n
  E[(\xi^{(n)}_k)^2] = 1+o(1)$ (from variance decomposition near (\ref{Sn_decomp}) and Lemma \ref{negligibility}).
\qed

\sect{Proof of Variance Lower Bound}
In this section, 
we prove Proposition 
\ref{lowerbound}.

\begin{lemma}
Let $f$ and $g$ be measurable functions on  $(\X, {\cal B}(\X))$.  Let
$\lambda$ be a probability  measure on $\X\times {\X}$ with marginals
$\alpha$ and $\beta$ respectively. Let $\pi (x_1, dx_2)$ and ${\widehat\pi}(x_2, dx_1)$ be the transition probabilities in the two directions
so that
$$
\alpha\pi=\beta,\quad \beta{\widehat\pi}=\alpha
$$
If 
$$\int f(x_1)\alpha(dx_1)=\int g(x_2)\beta(dx_2)=0$$
then,
$$
\bigg|\int f(x_1)g(x_2)\lambda(dx_1, dx_2)\bigg |\le\sqrt {\delta(\pi)} \,\|f\|_{L_2(\alpha)}\|g\|_{L_2(\beta)}
$$
\end{lemma}

{\it Proof.} Let us  construct a measure on ${\X}\times{\X}\times{\X}$
by starting with $\lambda$ on ${\X}\times{\X}$ and using reversed
$\widehat{\pi}(x_2, dx_3)$ to go from $x_2$ to $x_3$. The transition probability from $x_1$ to $x_3$
defined by
$$
Q(x_1,  A)=\int {\pi} (x_1, dx_2) \widehat{\pi}(x_2, A)
$$
satisfies $\delta(Q)\le \delta(\pi)$. Moreover $\alpha Q=\alpha$ and the
operator $Q$ is self adjoint and bounded with norm $1$ on
$L_2(\alpha)$. 
Then, if $f$ is a bounded function with
$\int f(x)\alpha(dx)=0$ (and so $E_\alpha[Q^n f]=0$), we have for $n\geq 1$,
\be
\label{bd}
\|Q^n f\|_{L_2(\alpha)}\le \|Q^n f\|_{L_\infty}\le (\delta(Q))^n {\rm
  Osc}(f).\ee
Hence, as bounded functions are dense, on the subspace
of functions, $M=\{f\in L_2(\alpha): \int f(x)\alpha(dx)=0\}$, the top
of the spectrum of $Q$ is less than $\delta(Q)$ and so
$\|Q\|_{L_2(\alpha,M)} \leq 
\delta(Q)$.  Indeed, suppose the
spectral radius of $Q$ on $M$ is larger than $\delta(Q)+\epsilon$ for
$\epsilon>0$, and $f\in M$ 
is a non-trivial bounded function 
whose spectral decomposition is with respect
to spectral values larger than $\delta(Q) + \epsilon$.
Then, $\|Q^nf\|_{L^2(\alpha)}\geq
\|f\|_{L^2(\alpha)}(\delta(Q)+\epsilon)^n$ which
contradicts the bound (\ref{bd}) when $n\uparrow \infty$.
[cf. Thm. 2.10 \cite{Seneta_book} for a proof in discrete
space settings.]
 
Then,
$$
\|\widehat{\pi} f\|^2_{L_2(\beta)}=<\pi\widehat{\pi} f,
f>_{L_2(\alpha)}=<Qf,f>_{L_2(\alpha)}\le \|Q\|_{L_2(\alpha,M}\|f\|_{L_2(\alpha)}
\le
\delta(Q)\|f\|^2_{L_2(\alpha)}.
$$
Finally 
$$
|\int f(x_1)g(x_2)\lambda(dx_1, dx_2)|=|<\widehat{\pi}f, g>_{L_2(\beta)}|\le 
\sqrt {\delta(\pi)} \,\|f\|_{L_2(\alpha)}\|g\|_{L_2(\beta)}.
$$

\qed
\begin{lemma}
\label{pairbound}
Let $f(x_1)$ and $g(x_2)$ be square integrable with respect to $\alpha$ and $\beta$ respectively.
Then,
$$E^\lambda \big[\big(f(x_1) - g(x_2)\big)^2\big] \ \geq \ \alpha(\pi)\
V\big(f(\cdot)\big)$$
as well as 
$$E^\lambda \big[\big(f(x_1) - g(x_2)\big)^2\big] \ \geq \ \alpha(\pi)\
V\big(g(\cdot)\big)$$

\end{lemma}

{\it Proof.}
We can assume without loss of generality that $f$ and $g$ have mean $0$ with respect to $\alpha$ and $\beta$ respectively. Then

\begin{align*}
E^\lambda \big[\big(f(x_1) - g(x_2)\big)^2\big]&=E^\alpha\big[ [f(x_1)]^2\big]+
E^\beta\big[ [g(x_2)]^2\big]-2 E^\lambda\big[ f(x_1)g(x_2)\big]\\
&\ge E^\alpha\big[ [f(x_1)]^2\big]+
E^\beta\big[ [g(x_2)]^2\big]-2 \sqrt {\delta(\pi)} \,\|f\|_{L_2(\alpha)}\|g\|_{L_2(\beta)}\\
&\ge (1-\delta(\pi))\|f\|^2_{L_2(\alpha)} \ =\ \alpha(\pi)\|f\|^2_{L_2(\alpha)}
\end{align*}
The proof of the second half is identical. \qed

{\vskip .1cm}
{\it Proof of Proposition \ref{lowerbound}.}
Applying Lemma \ref{pairbound} to the Markov pairs $\{(X^{(n)}_k,
X^{(n)}_{k+1}): 1\leq k\leq n-1\}$ with $f(X^{(n)}_k)=
E[Z^{(n)}_{k+1}|X^{(n)}_k]$ and $g(X^{(n)}_{k+1}) = Z^{(n)}_{k+1}$,
we get 
$$
E\big[\big(Z^{(n)}_{k+1} - E[Z^{(n)}_{k+1}|X^{(n)}_k]\big)^2\big]\ge \alpha_n E\big[\big(Z^{(n)}_{k+1}\big)^2\big] 
$$
On the other hand from (\ref{rearrange}), for $1\leq
k\leq n-1$, we have
\begin{align*}
V(f^{(n)}_k(X^{(n)}_k) &\leq E\big[
\big(f^{(n)}_k(X^{(n)}_k)\big)^2\big]\\
&\leq 2 E\big[ \big(Z^{(n)}_{k}\big)^2\big]+2E\big[\big(E\big[ Z^{(n)}_{k+1}|X^{(n)}_k\big]\big)^2\big]\\
&\leq 2 E\big[ \big(Z^{(n)}_{k}\big)^2\big]+2E\big[\big(Z^{(n)}_{k+1}\big)^2\big]
\end{align*}
Summing over $k$, and noting variance decomposition near (\ref{Sn_decomp}),
\begin{align*}
\sum_{k=1}^{n}V(f^{(n)}_k(X^{(n)}_k)
&\le 4\sum_{k=1}^n E\big[ \big(Z^{(n)}_{k}\big)^2\big]\\
&\le \frac{4}{\alpha_n}\bigg[\sum_{k=1}^{n-1} E\big[\big(Z^{(n)}_{k+1}
- E[Z^{(n)}_{k+1}|X^{(n)}_k]\big)^2\big] + E[(Z_1^{(n)})^2]\bigg]
 =  \frac{4}{\alpha_n}V(S_n).
\end{align*} 
 \qed

\bibliographystyle{plain}

\vskip .5cm

Sunder Sethuraman

400 Carver Hall

Department of Mathematics

Iowa State University

Ames, IA  50011

sethuram@iastate.edu

\vskip .5cm

S.R.S. Varadhan

251 Mercer St.

Courant Institute

New York University

New York, NY 10012

varadhan@cims.nyu.edu

\end{document}